
\documentstyle{amsppt}
\magnification=\magstep1
\pageheight{7in}
\TagsOnRight
\NoRunningHeads
\NoBlackBoxes

\def\tt#1{\medskip\leftline{\bf #1}\medskip\flushpar\hglue0pt}
\def\tit#1{\par\medskip\noindent{\bf #1}\par\smallskip\noindent\hglue0pt}

\def\ci#1{_{{}_{\ssize #1}}}
\def\cci#1{_{{}_{\sssize #1}}}
\def\ui#1{^{{}^{{}_ {\,#1} }}}

\def\R{\Bbb R}
\def\C{\Bbb C}
\def\D{\Cal D}

\def\T{\Bbb T}

\def\Vol{\operatorname{Vol}}
\def\Area{\operatorname{Area}}

\redefine\ge{\geqslant}
\redefine\le{\leqslant}

\def\la{\lambda}

\def\f{\varphi}
\def\e{\varepsilon}

\def\s{\sigma}   

\def\d{\delta}

\def\<{\langle}
\def\>{\rangle}

\topmatter

\title
Local dimension-free estimates for volumes of sublevel sets 
of analytic functions
\endtitle
\thanks The research was partially supported by the United States - Israel
Binational Science Foundation \endthanks

\author
F. Nazarov, M. Sodin, A. Volberg
\endauthor

\address (F.N.:) Department of Mathematics, Michigan State
University, East Lansing, MI 48824, U.S.A.
\endaddress
\email fedja\@math.msu.edu \endemail

\address (M.S.:) School of Mathematical Sciences, Tel Aviv University,
Ramat Aviv, 69978, Israel
\endaddress
\email sodin\@post.tau.ac.il \endemail  

\address (A.V.:) Department of Mathematics, Michigan State
University, East Lansing, MI 48824, U.S.A.
\endaddress
\email volberg\@math.msu.edu \endemail

\endtopmatter

\head
$\S 1$. The result
\endhead

\par\noindent 
In what follows, we denote complex balls $\{z\in \C^n:\, |z-w|<r\}$ 
by  $B_c(w,r)$ and real balls $\{x\in \R^n: |x-u|<r\}$
by $B(u,r)$.
For any real ball $B$, we denote by 
$\Vol\ci{B}$ the normalized volume
$\Vol\ci{B}(E) = \frac{\Vol (B\cap E)}{\Vol (B)}$.

We shall prove  

\tit{Theorem}{\it 
Let $F$ be a non-constant analytic function in $B_c(0,1)\subset\C^n$
such that $\sup_{B_c(0,1)}|F|\le~1$. Let $\e\le
\frac{1}{4}$, let $B$ be any real ball contained in $B(0, 1-\e)$, and
let $M\ci{B}(F)$ be a (unique) positive number  such that
$$
\Vol\ci{B}\{|F|\ge M\ci{B}(F)\} = \frac{1}{e}\,.
$$
Then, for every $\la> 1$,
$$
\Vol\ci{B}\{|F|\le (C\la)^{-\s} M\ci{B}(F)\} \le  \frac{1}{\la}\,,
\tag 1.1
$$
and
$$
\Vol\ci{B}\{|F|\ge (C\la)^\s M\ci{B}(F)\} \le e^{-\la}\,,
\tag 1.2
$$
where one can take $C=8$ and $\s=48\e^{-3}\,\log\frac{1}{|F(0)|}$.
}

\medskip

The main feature of the result is its dimensionless
character. Dimension-dependent versions of the theorem were obtained by
N.~Garofalo and P.~Garrett (see [GG]), and A.~Brudnyi (see [Br1], [Br2]).

If needed, the reader can adjust the theorem to
plurisubharmonic functions in the unit ball of $\C^n$ 
and to analytic functions with values in a Banach space.
Without changing the proof, one can replace real balls
$B$ by arbitrary convex bodies $V\subset B(0,1-\e)$
whose boundaries have sectional curvatures bounded from below
by some fixed positive constant. This ``curvature" restriction can,
probably, be relaxed but cannot be removed completely: a
simple example given in the end of this
note shows that estimates (1.1) and (1.2) may fail for thin
rectangles in $B(0,1-\e)\subset \R^2$.

Compiling the theorem with the technique from [NSV, \S3],
one can obtain an Offord-type statement about
the distribution of zeroes of analytic functions in families
that depend analytically on some parameters.
Informally speaking, the result is that the portion of the
family occupied by the functions
whose distribution of zeroes deviates from the ``average'' one by some
fixed amount, is about $\operatorname{Const}\,\exp\{-\operatorname{size\
of\ the\ deviation}\}$.
This might be a possible embryo of a non-linear and dimensionless
value-distribution theory. 

The theorem appeared as an attempt to ``generalize" the similar
statement for polynomials $P$ in $\R^n$. The main difference
is that, for polynomials, the counterparts of (1.1) and (1.2) hold with
$\s=\operatorname{deg} P$ and $C=4$ in {\it any} convex body $V\subset\R^n$
(see [NSV]). The quantity $\log \frac{1}{|F(0)|}$
appears as a ``natural analogue" of the degree of a polynomial just as it
does in the classical Cartan lemma.

As to the history of ``dimension-free estimates", the pioneering
dimensionless results are due to A.~C.~Offord [O], M.~Gromov and V.~Milman
[GM1] (the case of linear functions), and J.~Bourgain [B] (a
somewhat cruder form of (1.2) for polynomials). For
other developments, see A.~Brudnyi [Theorem~1.11, Br2]), S.~Bobkov
[Bo], and Carbery and Wright [CW].

\medskip
The proof of the theorem will be cooked from three ingredients. 

\tt{A. The geometric Kannan-Lov\'asz-Simonovits lemma:}
A continuous function $\Phi:\R^n\to \R_+$ is called
{\it logarithmically concave} if
$$
\Phi\left( \frac{x+y}{2}\right)\ge \sqrt{\Phi(x)\,\Phi(y)}
$$
for all $x,y \in \R^n$. 

\tit{Lemma~A} 
{\it Let $\Phi$ be a logarithmically concave function in $\R^n$.
Let $S\subset \text{supp}(\Phi)$ be a convex compact, and let $E\subset S$ be
a closed subset. For $\la >1$, define
$$
E_{\la, S}:= \Bigl\{
x\in E\,:\,
\frac{|E\cap J|}{|J|}\ge \frac{\la-1}\la
\text{ for every interval $J$ such that }x\in J\subset S,\
\Bigr\}.
$$
Then
$$
\frac{\int_{E_{\la, S}}\Phi }{\int_S \Phi}
\le
\left(
\frac{\int_E \Phi}{\int_S \Phi}\right
)^\la\,.
$$
}

This lemma was proved in [NSV] using the needle
decomposition technique developed by M.~Gromov and V.~Milman [GM2] and by
L.~Lov\'asz and M.~Simonovits [LS]. It can also be derived
from a result of
R.~Kannan, L.~Lov\'asz, and M.~Simonovits [Theorem~2.7, KLS]. 

\tit{B. One dimensional Remez property:}
We shall use the following result (which, probably, should
be called the Boutroux-Cartan-Remez property): 
\tt{Lemma~B}
{\it
Let $f$ be an analytic function in the unit disk $\Bbb D$ such that 
$\sup_{\Bbb D}|f| \le 1$, and let $|f(a)|=|f(-a)|>0$ for some $a\in(0,1)$.
Then $f$ has the Remez property on the interval $[-a,a]$, i.e.,
for every sub-interval $I\subset [-a,a]$ and every set $E\subset I$,
$$
\max_I |f| \le 
\left( 
\frac{C|I|}{|E|}\right)^\s \sup_E |f|
$$
with $C=8$ and $\s=\frac{3}{1-a}\, \log\frac{1}{|f(a)|}$.
}

\medskip\par\noindent 
For the sake of
completeness, we provide the proof of Lemma~B in $\S 2$.

\tt{C. A change of variable:}
Let
$\d\le \frac{1}{8}$.
Set
$A=1-\d^3$, $a=\sqrt A$,  
$\f(\zeta)=\frac{A-\zeta}{1-A\zeta}$,
and
consider the mapping $T$ defined on the unit ball $B_c(0,1)$ in $\C^n$ by 
the formula
$$\tsize
T(z):=\f\Bigl(\sum\limits_{j=1}^n z_j^2\Bigr)z\,,\qquad z=(z_1,\dots,z_n)\in
B_c(0,1)\subset \C^n\,.
$$
In particular, $T(x)=\f(|x|^2)x$ for $x\in B(0,1)$. 

\tt{Lemma~C}
{\it Set
$R_0=1-3\d-\d^3$,
$r_0=\sqrt R_0$. Then 
the mapping $T$ has the following properties:
\roster
\item
$TB_c(0,1)\subset B_c(0,1)$;
\item
$T$ maps the real sphere $|x|=a$ to the origin;
\item
$T$ is one-to-one in the ball $B(0,r\ci 0)$;
\item
$TB(0,r\ci0)$ is a ball centered at the origin of radius greater than
$1-2\d$;
\item
The Jacobian $|\det D_x T|$ is a logarithmically concave function
in $B(0,r\ci0)$;
\item  
The (partial) pre-image $B(0,r\ci 0)\cap T^{-1}B$
of every (real) ball $B\subset TB(0,r\ci0)$ is convex.
\endroster
}
\medskip\par\noindent The first two properties are obvious; the others 
will be proved in $\S 3$.

\tt{Proof of Theorem:}
Let $F$ be a non-constant analytic function in $B_c(0,1)$.
We shall show that for every $c>0$ and every $\la >1$,
$$
\Vol\ci{B} \{|F|\ge (C\la)^\s c\}
\le \Big( \Vol\ci{B} \{|F|\ge c \} \Big)^\la\,.
\tag 1.3
$$
The rest is the same as in [NSV]: to get (1.2), we
just set $c=M\ci{B}(F)$ in (1.3);
to get (1.1), we rewrite (1.3) in the form
$$
\Vol\ci{B} \{|F|\ge c \} \le \Big( 1 - \Vol\ci{B} \{|F| < (C\la)^{-\s} c \}
\Big)^\la
$$ 
and, taking $c=M\ci{B}(F)$, obtain
$$
\Vol\ci{B} \{|F| < (C\la)^{-\s} M\ci{B}(F)\} \le 1 - e^{-1/\la} <
\frac{1}{\la}\,,
$$
which is identical to (1.1) since $\Vol\ci{B}\{|F|=\text{const}\}=0$.

To prove (1.3), choose $\d=\frac{\e}{2}$ and
consider the composition $F\ci{T}(z)=(F\circ T)(z)$
of the function $F$ with the mapping $T$
defined above. The function $F\ci{T}$ is analytic in the complex
unit ball and $\sup_{B_c(0,1)} |F\ci{T}|\le 1$.
The advantage we gain by this trick is that the new function
$F\ci{T}$ has a lower bound on a
massive set (the real sphere) instead of just one point (the origin):
$F\ci{T}(u)=F(0)$ for every $u\in \R^n$ with $|u|=a$.
Let $S=B(0,r\ci 0)\cap T^{-1}B$.
Due to Lemma~C (property (6)), this is a convex compact
subset of $B(0,r_0)$. 
We shall show that for every $c>0$ and for every $\la >1$,
$$
\frac{\int_{S\cap\{|F\cci{T}|\ge (C\la)^{{}^\s} c\}}
|\det D_x T|}
{\int_S |\det D_xT|}
\le
\left(
\frac{\int_{S\cap\{|F\cci{T}| \ge c\}} |\det D_x T|}
{\int_S |\det D_xT|} \right)^\la
\tag 1.4
$$
which is equivalent to (1.3).

Let $E=\{x\in S:\, |F\ci{T}(x)|\ge c\}$.
To prove (1.4), we check that 
$$
S\cap \{|F\ci{T}|> (C\la)^\s c\} \subset
E_{\la, S}\,,
\tag 1.5
$$
where the set $E_{\la, S}$ is defined in Lemma~A. 
Since $F\ci{T}$ is a non-constant analytic function, the level set
$\{|F\ci{T}|=(C\la)^\s c\}$ has zero volume. Then Lemma~A with
the function $\Phi = |\det D_xT|$ (which is logarithmically
concave due to property~(5) in Lemma~C) gives us (1.4).

Assume that $x\notin E_{\la, S}$, i.e., that
there exists an interval $J\subset S$ containing the point $x$ and such
that the length of the set $J\setminus E$ is at least 
$\la^{-1} |J|$. Extend this interval until the endpoints appear
on the unit sphere $\partial B(0,1)$ and denote the extended
interval by $J^*$.
Let $\Delta$ be the one-dimensional complex disk with
diameter $J^*$.
Then $\Delta \subset B_c(0,1)$ and
$|F\ci{T}(x)|=|F(0)|$ for $x\in J^*\cap \partial B(0,a)$. Further,
$|J^*\cap B(0,a)|\le a|J^*|$ and we can apply the 
one-dimensional Remez property (Lemma~B) to the analytic function 
$F\ci{T}\big|_\Delta$, the interval $J$, and its subset $J\setminus
E$. We get
$$
|F\ci{T}(x)| \le \max_J|F\ci{T}| \le
\left( \frac{C|J|}{|J\setminus E|} \right)^\s 
\sup_{J\setminus E} |F\ci{T}|
\le (C\la)^\s c\,,
$$
with $C=8$ and 
$$
\split
\s = \frac{3}{1-a}\log\frac{1}{|F\ci{T}(a)|} =
\frac{3}{1-\sqrt{1-\d^3}} \log\frac{1}{|F(0)|}
\\
< \frac{6}{\d^3} \log\frac{1}{|F(0)|} =
\frac{48}{\e^3}\log\frac{1}{|F(0)|}\,,
\endsplit
$$ 
completing the proof of (1.5) and, thereby, of the theorem.
$\square$

\head
$\S 2$. Proof of Lemma~B
\endhead
We shall use the standard factorization $f(z)=U(z)B(z)$ where $U(z)$ has 
no zeroes in the disk and $B(z)$ is the Blaschke product.
Since for every $x\in[-a,a]$,
$$
\log |U(x)|=-\int_{\T}\frac{1-x^2}{|1-x\zeta|^2}d\mu(\zeta)\,
$$
where $\mu$ is some positive measure on the unit circle $\T$, and
since
$$
\frac{1}{|1-x\zeta|^{2}}\le
\frac{1}{|1-a\zeta|^2}+
\frac{1}{|1+a\zeta|^2}
$$
for every $\zeta\in\D$, $x\in[-a,a]$,
we immediately conclude that
$$
\log |U(x)| \ge -\frac{1-x^2}{1-a^2} \int_{\T} 
\left( \frac{1-a^2}{|1-a\zeta|^2} + \frac{1-a^2}{|1-a\zeta|^2}
\right) d\mu(\zeta) = \frac{1-x^2}{1-a^2} \log|U(a)U(-a)|\,
$$
and, therefore,
$$
\min_{[-a,a]}|U|\ge 
\bigl|U(-a)U(a)\bigr|\ui{\tfrac{1}{1-a^2}}\,.
\tag 2.1
$$
We shall split the Blaschke product $\dsize
B(z)=\prod_{\zeta}\frac{z-\zeta}
{1-z \bar \zeta }$ into two factors: $B\ci1(z)$, which is the product over 
all
zeroes $\zeta$ satisfying
$$
\frac{1-|\zeta|^2}{|1+a\zeta|^2}+
\frac{1-|\zeta|^2}{|1-a\zeta|^2}\le \frac{2}{3}\,,
$$
and $B\ci2(z)$, which is the product over all zeroes $\zeta$ for which
the opposite  inequality holds.
Our next aim will be to show that for all $x\in[-a,a]$,
$$
|B\ci1(x)|\ge\bigl|B\ci1(-a)B\ci1(a)\bigr|\ui{\frac{2(1-x^2)}{1-a^2}}\,,
$$
which yields
$$
\min_{[-a,a]}|B_1| \ge \left|
B_1(a)B_1(-a)
\right|^{\frac{2}{1-a^2}}\,.~
\tag 2.2 
$$
Clearly, it is enough to establish this inequality for every Blaschke 
factor
in $B\ci1(z)$. Using the inequality $1-t\ge e^{-2t}$ ($0\le t\le
\frac{2}{3}$),
we obtain
$$
\split
\left|\frac{x-\zeta}{1-x\bar\zeta}\right|^2 &=
1-\frac{(1-x^2)(1-|\zeta|^2)}{|1-x\bar\zeta|^2} \\
&\ge
1 - \left[
\frac{(1-x^2)(1-|\zeta|^2)}{|1+a\bar\zeta|^2}+
\frac{(1-x^2)(1-|\zeta|^2)}{|1-a\bar\zeta|^2}
\right] \\ 
&\ge
\exp
\left\{
-\frac{2(1-x^2)}{1-a^2} 
\left[ 
\frac{(1-a^2)(1-|\zeta|^2)}{|1+a\bar\zeta|^2}
+ \frac{(1-a^2)(1-|\zeta|^2)}{|1-a\bar\zeta|^2}
\right] \right\} \\  
&\ge
\left[1-\frac{(1-a^2)(1-|\zeta|^2)}{|1+a\bar\zeta|^2}\right]
^{\tfrac{2(1-x^2)}{1-a^2}}\cdot
\left[1-\frac{(1-a^2)(1-|\zeta|^2)}{|1-a\bar\zeta|^2}\right]
^{\tfrac{2(1-x^2)}{1-a^2}} \\
&=
\left|\frac{-a-\zeta}{1+a\bar\zeta}
\cdot   
\frac{a-\zeta}{1-a\bar\zeta}\right|
^{\frac{2(1-x^2)}{1-a^2}},
\endsplit
$$
proving the statement.

The next observation is that the number $N$ of factors in $B\ci2(z)$
satisfies the inequality
$$
N \le \frac3{1-a^2} \log\frac{1}{|B\ci2(-a)B\ci2(a)|}\,.
\tag 2.3
$$
Indeed, for every zero $\zeta$ in $B\ci 2$, we have
$$      
\split
\left|\frac{-a-\zeta}{1+a\bar\zeta}
\cdot
\frac{a-\zeta}{1-a\bar\zeta}\right|^2 &=
\left[1-\frac{(1-a^2)(1-|\zeta|^2)}{|1+a\bar\zeta|^2}\right]\cdot
\left[1-\frac{(1-a^2)(1-|\zeta|^2)}{|1-a\bar\zeta|^2}\right] \\
&\le\exp\left\{-(1-a^2)
\left[\frac{1-|\zeta|^2}{|1+a\bar\zeta|^2}+
\frac{1-|\zeta|^2}{|1-a\bar\zeta|^2}\right]\right\} \\
&\le e^{\tsize -\frac{2(1-a^2)}{3}}\,.
\endsplit
$$   
Thus,
$$
|B_2(a)B_2(-a)| \le \exp\left[-\frac{N(1-a^2)}{3} \right]\,,
$$
which is equivalent to (2.3).

Now write $B_2(z)=P(z)R(z)$ where $P(z)=\prod_{k=1}^N (z-\zeta_k)$,
and $R(z) = \prod_{k=1}^N \dfrac{1}{1-z\bar\zeta_k}$.
We have
$$
\max_{[-a.a]} |R| 
\le \left( \frac{1+a}{1-a} \right)^N \min_{[-a,a]} |R| 
\le \left( \frac{2}{1-a} \right)^N \min_{[-a,a]} |R|\,.
\tag 2.4
$$
At last, according to the classical Remez inequality (see, for example,
[DR] or [BG]), for any sub-interval 
$I\subset [-a,a]$ and any measurable subset $E\subset I$,
$$
\max_I |P| 
\le \left( \frac{4|I|}{|E|}
\right)^N \sup_E |P|\,.
\tag 2.5
$$
Combining estimates (2.1)--(2.5), we get
$$
\split
\max_I |f| &\le \max_I |B_2| \\
&\le \max_I |P| \cdot \max_{[-a,a]} |R| \\
&\le 
\left(
\frac{4|I|}{|E|}
\right)^N 
\sup_E |P| 
\cdot 
\left(\frac{2}{1-a} \right)^N 
\min_{[-a,a]} |R| \\ 
&\le 
\left(
\frac{8|I|}{|E|}
\right)^N 
\cdot 
\left(\frac{1}{1-a} \right)^N
\sup_E |B_2| \\
&\le 
\left(\frac{8|I|}{|E|}\right)^N 
\cdot \left(\frac{1}{1-a} \right)^N 
\cdot \max_{[-a,a]} \frac{1}{|U|} 
\cdot \max_{[-a,a]} \frac{1}{|B_1|}
\cdot \sup_E |f| \\
&\le 
\left(\frac{8|I|}{|E|}\right)^N 
\cdot \left(\frac{1}{1-a} \right)^N 
\cdot
\left| \frac{1}{U(a)U(-a)B_1(a)B_1(-a)}
\right|^{\frac{2}{1-a^2}}
\sup_E |f| \\
&\le 
\left(\frac{8|I|}{|E|}\right)^\s 
\sup_E |f|\,,
\endsplit
$$
where
$$
\split
\s &\le 
\left( 1 + \log\frac{1}{1-a} \right) N
+ \frac{2}{1-a^2} \log\frac{1}{|U(a)U(-a)B_1(a)B_1(-a)|} \\
&\le 3 
\left( 1 + \log\frac{1}{1-a} \right) 
\log\frac{1}{|B_2(a)B_2(-a)|}
+ \frac{2}{1-a} \log\frac{1}{|U(a)U(-a)B_1(a)B_1(-a)|} \\
&\le\frac{3}{1-a} \log\frac{1}{|f(a)f(-a)|}\,
\endsplit
$$
(in the last line we used the inequality $1+\log t\le t$ when $t\ge 1$).
Lemma~B is proved.
$\square$

\head
$\S 3$. Proof of Lemma~C
\endhead
\tt{$T$ is one-to-one in the ball $B(0,r_0)$:}
We show that the function $r\mapsto r\f(r^2)$ where, as before,
$\f (\zeta) = \frac{A-\zeta}{1-A\zeta}$,
is increasing on the interval $[0,r_0]$. Set $R=r^2$.  
We have
$$
\frac{d}{dr}(r\f(r^2))=\f(R)\left(1+2R\frac{\f'(R)}{\f(R)}\right).
$$
Since $0\le R\le R\ci0<A$, we have $\f(R)>0$. So, it will suffice to show 
that $\frac{|\f'(R)|}{\f(R)}\le \frac{1}{2}$. A direct computation yields
$$
\frac{|\f'(R)|}{\f(R)}\le
\frac{|\f'(R)|}{\f(R)^2}= \frac{1-A^2}{(A-R)^2}\le
\frac{2(1-A)}{(A-R_0)^2}\le \frac{2\d}{9} < \frac{1}{30}\,,
$$
since $\d\le\frac{1}{8}$. $\square$

\tt{$TB(0,r_0)$ is a ball centered at the origin 
with radius bigger than $1-2\d$:}
It is clear now that $TB(0,r\ci0)=B(0,r\ci0\f(R\ci0)\,)$, so we need only 
to show that $r\ci0\f(R\ci0)>1-2\d$.
We have
$$
1-AR\ci0=1-(1-\d^3)(1-3\d-\d^3)= 3\d +\d^3+\d^3(1-3\d-\d^3)\le
3\d+2\d^3
$$
and, thereby,
$$
\f(R\ci0)=
\frac{A-R\ci0}{1-AR\ci0}\ge 
\frac{3\d}{3\d+2\d^3} > \frac{1}{1+ \d^2}\,.
$$
Thus, to prove our inequality, we need to check that
$$
1-3\d-\d^3\ge (1-2\d)^2\bigl(1+\d^2\bigr)^2.
$$
The right hand side does not exceed
$$
(1-4\d+4\d^2)(1+3\d^2) \le 1-4\d + 7\d^2\,.
$$
Since $\d\le \frac{1}{8}$, we have $7\d^2 + \d^3 < 8\d^2 \le \d$, 
 finishing the proof. $\square$

\tt{The Jacobian $|\det D_xT|$ is a
logarithmically concave function in $B(0,r_0)$:}
First, we compute the Jacobian. Let $T_i(x)=\f(|x|^2)x_i$. 
Then
$$
\frac{\partial T_i}{\partial x_j} = \cases 
\f'(r^2) 2x_ix_j, &i\ne j \\ \\
\f(r^2) + \f'(r^2)2x_i^2, &i=j
\endcases 
$$
whence $\det D_x T = \det(\xi I+A)$, where $\xi = \f(r^2)$ and 
$A_{i\,j}=2\f' (r^2)x_ix_j$. Since the rank of $A$ is one,
$\det(\xi I + A) = \xi^n + \xi^{n-1}\text{tr}(A)$,
and 
$$
|\det D_xT|=\big(\f(r^2)+2r^2\f'(r^2)\big)\cdot \big(\f(r^2)\big)^{n-1}.
$$
(This result can also be obtained in a purely geometric way:
just consider the image of a small domain containing $x$ and
bounded by two concentric spheres and
a thin cone).

The Taylor expansions
$$
\f(R)=A-(1-A^2)\sum_{k=1}^\infty A^{k-1}R^k\,,
$$
and
$$
\f(R) + 2R\f'(R) = A-(1-A^2) \sum_{k=1}^\infty (1+2k)A^{k-1}R^k
$$
immediately show that both 
$\f(r^2)$ and $\f(r^2)+2r^2\f'(r^2)$  
are concave decreasing functions of 
$r$ on the interval $[0, 1]$. Since they are also positive on 
$[0,r\ci 0]$,
they are logarithmically concave on that interval.
Hence the function $r\mapsto [\f(r^2)+2r^2\f'(r^2)]\,[\f(r^2)]^{n-1}$
is also logarithmically concave on the interval $[0,r_0]$.
It remains to recall that if $\Phi(r)$ is a decreasing logarithmically
concave function on the interval $[0,r\ci0]$, then $x\mapsto \Phi(|x|)$ is
logarithmically concave in the ball $B(0,r\ci 0)\subset\R^n$. $\square$

\tt{The pre-image $T^{-1}B$ of every (real) ball 
$B\subset TB(0,r_0)$ is convex:}
Since the pre-image $T^{-1}B$ is a body
of revolution around the axis containing both the origin and the center of
the ball $B$, it is enough to prove our statement on the plane $\R^2$.
In order to do so, we shall show that the curvature of the   
image of any straight line tangent to the boundary of $T^{-1}B$ does not
exceed the curvature of the boundary of $B$ which is 
$\frac{1}{\operatorname{rad}(B)}$. It is going to be a
simple but somewhat boring exercise in
differential geometry.
  
Let $rx$ ($0\le r\le r\ci0$, $x\in\R^2$, $|x|=1$) be a point on the 
boundary of $T^{-1}B$ and let $y(t)=rx+tv$ ($v\in \R^2$, 
$|v|=1$, $t\in\R$) be 
the corresponding tangent line. Let $\alpha$ be the angle between the 
vectors $x$ and $v$. The image of our tangent line is the curve
$$
\s(t)=\f(|y(t)|^2)y(t)=\f(r^2+2rt\cos\alpha +t^2)(rx+tv).
$$
To estimate the curvature, we need to compute the first and second
derivatives of $\s$. Differentiation yields
$$
\align
\s'(t)&=\f(|y(t)|^2)v+2\f'(|y(t)|^2)\langle y(t),v\rangle y(t)\,,
\\
\s''(t)&=4\f'(|y(t)|^2)\langle y(t),v\rangle v + 2\f'(|y(t)|^2) y(t)
+4\f''(|y(t)|^2)\langle y(t),v\rangle^2 y(t).
\endalign
$$
Plugging in $t=0$ and denoting, as above, $r^2=R$, we obtain
$$
\aligned
\s'(0)&=\f(R)v+2R\f'(R)(\cos\alpha)\, x\,,
\\
\s''(0)&=4r\f'(R)(\cos\alpha)\, v + 2 r\f'(R)x
+4rR\f''(R)(\cos^2\alpha)\, x.
\endaligned
$$
Now we are ready to estimate the curvature. We shall use the
standard formula
$$
\operatorname{curvature}=\frac{|\s'(0)\times\s''(0)|}{|\s'(0)|^3}\,.
$$
We have
$$
\frac{|\s'(0)|}{\f (R)}
\ge 1 - 2R \frac{|\f'(R)|}{\f (R)}\cos\alpha
\ge \frac{14}{15} 
$$
(recall that $R<1$ and $\frac{|\f'(R)|}{\f(R)}\le\frac{1}{30}$),
and therefore
$|\s'(0)|^{3}\ge \frac{4}{5} \f(R)^3$.
Using this estimates, we finally obtain
$$
\split
\operatorname{curvature} &\le
\frac{|\s'(0)\times\s''(0)|}{\frac{4}{5}\f(R)^3} \\ 
&= \frac{5}{2} r|v\times
x|\cdot\left|
\frac{\f'(R)}{\f(R)^2}+2R\left[
\frac{\f''(R)}{\f(R)^2}-\frac{2\f'(R)^2}{\f(R)^3}
\right] \cos^2\alpha 
\right| 
\\
&= \frac{5}{2} r|\sin\alpha|\cdot\left|
\frac{\f'(R)}{\f(R)^2}+2R\left[
\frac{\f'(R)}{\f(R)^2}
\right]' \cos^2\alpha 
\right| \\
&\le \frac{5}{2} \left[
\frac{1-A^2}{(A-R)^2}+4R \frac{1-A^2}{(A-R)^3}
\right] \\
&\le \frac{10(1-A^2)}{(A-R)^3}\cdot\left[1+\frac{A-R}{4}
\right] \\
&\le \frac{20(1-A)}{(A-R)^3}\cdot \frac{5}{4} \\
&= \frac{25}{27} < 1 < \frac{1}{\operatorname{rad}(B)}\,, 
\endsplit
$$
completing the proof of Lemma~C. $\square$

\head
$\S 4$. An example
\endhead

Let $Q(z)$ be an arbitrary polynomial. Let $\eta>0$ be so small that 
$$\eta \max_{|z|\le 1}|Q(z)|<\frac{1}{8}\,.$$
Consider the analytic function $F$ in the unit ball $B_c(0,1)\subset
\C^2$ defined by
$$
F(z_1, z_2) = \frac12\left[2\eta\, Q(z_1) + z_2 + \frac{1}{2}\right]
$$
and take rectangles
$$
V_\delta = \left\{0\le x_1\le \frac{1}{4}\,, \quad 0\le x_2+\frac{1}{2} \le \d
\right\} \subset B\left(0,\tfrac{3}{4}\right)\,,
\qquad 0<\d\le \frac{1}{2}\,.
$$
It is easy to see that $|F|\le 1$ in $B_c(0,1)$ and $|F(0,0)|\ge \frac14$
regardless of 
the choice of $Q$.  
Notice that for very small $\delta>0$, the distribution of $|F|$ in the
rectangle $V_\delta$
with respect to the normalized area $\frac{1}{\Area(V)}\,d\!\Area(x)$
is practically indistinguishable from the distribution of $\eta\, Q(t)$ on
the interval 
$[0,\frac14]$ with respect to the normalized Lebesgue measure $4dt$.
If the estimates (1.1) and (1.2) of the theorem were true in every
rectangle $V_\delta$, they would 
also hold for the measures of level sets of the polynomial $\eta\, Q(t)$ on
the interval $[0,\frac14]$. 
Since they
are scale-invariant, they would also hold for the measures of level sets of
the polynomial $Q$ on the 
interval $[0,\frac14]$. But, since polynomials are dense in the space of
continuous functions, this
would imply that they hold for level sets of {\it any} continuous function
$g(t)$ on the interval
$[0,\frac14]$, which is clearly false.
$\square$

\bigskip
\centerline{\bf References}
\bigskip

\medskip\par\noindent{\bf [Bo]} S. G. Bobkov,
{\it Remarks on the growth of 
$L^p$ norms of polynomials, }
Lecture Notes in Math. {\bf 1745} (2000), 27--35.

\medskip\par\noindent{\bf [B]} J. Bourgain,
{\it On the distribution of polynomials on high dimensional
convex sets, }
Lecture Notes in Math. {\bf 1469} (1991), 127--137.

\medskip\par\noindent {\bf [Br1]} A. Brudnyi, {\it Local
inequalities
for 
plurisubharmonic functions,} Ann. Math. {\bf 149} (1999), 511--533.

\medskip\par\noindent {\bf [Br2]} A. Brudnyi, {\it On the local
behavior
of analytic functions,} Journ. Funct. Anal. {\bf 169} (1999), 481--493.

\medskip\par\noindent {\bf [BG]} Yu. Brudnyi and M. Ganzburg,
{\it One extremal problem for polynomials of $n$ variables, } 
Izv. Akad. Nauk SSSR (Mat) {\bf 37} (1973), 344--355.
(Russian)

\medskip\par\noindent {\bf [CW]} A. Carbery and J. Wright,
{\it Distributional and $L^q$ norm
inequalities for polynomials over convex bodies in ${\Bbb R}^n$},
Math. Res. Lett. {\bf 8} (2001), 233--248.

\medskip\par\noindent {\bf [DR]} R. M. Dudley and B. Randol,
{\it Implications of pointwise bounds on polynomials,}
Duke Math. J. {\bf 29} (1962), 455--458.

\medskip\par\noindent {\bf [GG]} N.~Garofalo and P.~B.~Garrett, 
{\it $A\sb p$-weight properties
of real analytic functions in ${R}\sp n$,} Proc. Amer. Math. Soc. 
{\bf 96} (1986), 636--642.

\medskip\par\noindent{\bf [GM1]} M. Gromov and V. Milman,
{\it Brunn theorem and a concentration of volume phenomena for symmetric 
convex bodies,} Israel seminar on geometrical aspects of functional
analysis (1983/84), Tel Aviv University, Tel Aviv, 1984.

\medskip\par\noindent {\bf [GM2]} M. Gromov and V. Milman,
{\it Generalization of the spherical
isoperimetric inequality
to uniformly convex Banach spaces,} Compositio
Math. {\bf 62} (1987), 263-282.

\medskip\par\noindent {\bf[ KLS]} R. Kannan, L. Lov\'asz and
M. Simonovits,
{\it Isoperimetric problem for convex bodies and a localization lemma, }
Discrete and Comput. Geometry {\bf 13} (1995), 541--559.

\medskip\par\noindent {\bf [LS]} L. Lov\'asz and  M. Simonovits,
{\it Random walks in a 
convex body and an improved volume algorithm,} Random Structures and
Algorithms {\bf 4} (1993), 359--412.

\medskip\par\noindent {\bf [NSV]} F. Nazarov, M. Sodin and 
A. Volberg, {\it The geometric Kannan-Lov\'asz-Simonovits lemma, 
dimension-free estimates for volumes of sublevel sets of polynomials, 
and distribution of zeroes of random analytic functions,} Algebra \&
Analysis (St. Petersburg Math. Journ.), to appear.

\medskip\par\noindent {\bf [O]} A. C. Offord, {\it The distribution of
zeros of power series whose coefficients are independent random
variables}, Indian J. Math. {\bf 9} (1967), 175--196.

\enddocument